\documentclass[12pt,centertags,oneside]{amsart}

\usepackage{amsmath,amstext,amsthm,amscd,typearea,hyperref}
\usepackage{amssymb}
\usepackage{a4wide}
\usepackage[mathscr]{eucal}
\usepackage{mathrsfs}
\usepackage{typearea}
\usepackage{charter}
\usepackage{pdfsync}
\usepackage[a4paper,width=16.2cm,top=3cm,bottom=3cm]{geometry}

\numberwithin{equation}{section}

\newtheorem{theorem}{Theorem}[section]
\newtheorem{definition}[theorem]{Definition}
\newtheorem{proposition}[theorem]{Proposition}
\newtheorem{corollary}[theorem]{Corollary}

\newtheorem{problem}{Problem}
\newtheorem{conjecture}[problem]{Conjecture}

\newenvironment{problem1a}{\par\medskip
\noindent 
\textbf{Problem 1a.} \itshape
\noindent} {\medskip}

\newenvironment{problem1b}{ \par\medskip
\noindent 
\textbf{Problem 1b.} \itshape
\noindent}  {\medskip}

\newcommand{\cali}[1]{\mathscr{#1}}

\newcommand{\NS}{{\rm NS}}

\newcommand{\supp}{{\rm supp}}

\newcommand{\ddc}{{dd^c}}
\newcommand{\dc}{{d^c}}
\newcommand{\dbar}{{\overline\partial}}

\newcommand{\codim}{{\rm codim\ \!}}

\newcommand{\Cc}{\cali{C}}

\newcommand{\Ec}{\cali{E}}

\newcommand{\Hc}{\cali{H}}
\newcommand{\Kc}{\cali{K}}

\newcommand{\Uc}{\cali{U}}

\newcommand{\C}{\mathbb{C}}

\newcommand{\R}{\mathbb{R}}
\newcommand{\T}{\mathbb{T}}

\renewcommand{\P}{\mathbb{P}}

\title[Equidistribution problems of complex dynamics in higher dimension]{Equidistribution problems of complex dynamics\break in higher dimension}

\author{Tien-Cuong Dinh}
\address{Department of Mathematics, National University 
of Singapore, 10 Lower Kent Ridge Road, Singapore 119076.}
\email{matdtc@nus.edu.sg}
\thanks{T.-C. D.  supported by Start-Up 
Grant R-146-000-204-133 from National University of\break Singapore}

\author{Nessim Sibony}
\address{Laboratoire de Math\'ematiques d'Orsay, Univ. Paris-Sud, CNRS, Universit\'e
Paris-Saclay, 91405 Orsay, France 
and 
Korea Institute For Advanced Studies, 
Seoul, 130-722 South Korea.}

\email{Nessim.Sibony@math.u-psud.fr}

\date{November 07, 2016}

\begin{document}

\maketitle

\begin{center}
{\it to Ngaiming Mok for his sixtieth birthday}
\end{center}

\begin{abstract}
Equidistribution of the orbits of points, subvarieties or of periodic points in complex dynamics is a fundamental problem. It is often related to  strong ergodic properties of the dynamical system and to a deep understanding of analytic cycles, or more generally positive closed currents, of arbitrary dimension and degree. The later topic includes the study of the potentials and super-potentials of positive closed currents, their intersection with or without dimension excess. In this paper, we will survey some results and tools developed during the last two decades. Related concepts,  new techniques and open problems will be presented. 
\end{abstract}

\medskip\medskip

\noindent
{\bf Classification AMS 2010:} 32H, 32U, 37F.

 \medskip


\tableofcontents

\section{Introduction} \label{introduction}

Let $X$ be a compact K\"ahler manifold of dimension $k$ with a K\"ahler form $\omega$ on $X$. We will use the Riemannian metric on $X$ induced by this K\"ahler form. Let $f$ be a dominant meromorphic self-map or self-correspondence on $X$. Their formal definitions will be given in the next section. It is useful to notice that in general, $f$ has an indeterminacy set and then it may not be continuous. By correspondence, we mean a multi-valued map. 

The iterate of order $n$ of $f$ is roughly given by $f^n:=f\circ \cdots \circ f$, $n$ times. The adjoint correspondence is denoted by $f^{-1}$. Its graph in $X\times X$ is obtained as the image of the graph of $f$ by the involution $(x_1,x_2)\mapsto (x_2,x_1)$. In general, correspondences do not satisfy that $f\circ f^{-1}$ is the identity map on $X$.  Define $f^{-n}$ as the iterate of order $n$ of $f^{-1}$ which is also the adjoint correspondence of $f^n$. Note that the discussion is already interesting when $f$ is a meromorphic or holomorphic map but even in this case, $f^{-1}$ may not be a meromorphic map but  a meromorphic correspondence. So for convenience, we will work in the correspondence setting.

Periodic points of period $n$ are given by the intersection of the graph $\Gamma_n$ of $f^n$ with the diagonal $\Delta$ of $X\times X$, when we identify $\Delta$ with $X$ in the canonical way. This is an analytic subset of $X$. Since $\dim \Gamma_n=\dim \Delta=k$, at least in generic situations, we expect that $\Gamma_n\cap \Delta$ is of dimension 0, i.e., is a finite set. However, in general, this intersection may have positive dimension. Let $Q_n$ denote the set of isolated periodic points of period $n$, i.e., isolated points in the intersection $\Gamma_n\cap \Delta$. We will count the points in $Q_n$ with their multiplicities.
Periodic points are fundamental objects in dynamics. Since their orbits are simple (finite), they help us to understand more complicated orbits. Moreover, they are strongly related to the topological entropy of the map or correspondence and to the dependence of the system under a perturbation on the map, see \cite{BBD, DS10} and the references therein.
Here is the first basic question on periodic points.

\begin{problem}
Estimate the cardinality $\#Q_n$ of $Q_n$ when $n$ goes to infinity.
\end{problem}

When the set of periodic points of period $n$ of $f$ is finite, all of them are isolated and $\#Q_n$ can be computed using the classical Lefschetz fixed point formula. Equivalently, this number is equal to the intersection number of the cohomology classes $\{\Gamma_n\}$ and $\{\Delta\}$ associated to $\Gamma_n$ and $\Delta$ in $X\times X$. The computation is then easier in this particular case. When $\Gamma_n\cap \Delta$ is of positive dimension, it is not clear how to evaluate the contribution of the components of positive dimension in the last intersection.
We will see later the reason to consider the following question which is a particular case of the last problem.

\begin{problem1a} \it
Find a good upper bound for $\#Q_n$.
\end{problem1a}

More precisely, denote by $L(f^n)$ the Lefschetz number for $f^n$ which is the intersection number of the cohomology classes 
$\{\Gamma_n\}$ and $\{\Delta\}$ in $X\times X$. As mentioned above, this is the number of periodic points of period $n$ counting with multiplicity, in the case where these points are all isolated. So in this case, we have 
$$\#Q_n=L(f^n):=\{\Gamma_n\}\smallsmile \{\Delta\},$$
where the last cup-product is the intersection number of the cohomology classes associated to $\Gamma_n$ and $\Delta$.
We don't know if the following property holds in general.

\begin{problem1b} \it
Do we always have 
$$\#Q_n\leq L(f^n)+o(L(f^n)) \ \text{ as} \ n \text{ goes to infinity ?}$$
Or, find sufficient conditions on $f$ so that this property holds.
\end{problem1b}

A final goal in the dynamical study is to show that, under suitable hypotheses, the points in $Q_n$ are distributed according to a canonical probability measure when $n$ tends to infinity. One also has to understand the long 
term behavior of the dynamical system with respect to this probability measure. This is an analogue to 
the distribution of torsion points on 
Abelian varieties or other special points in algebraic geometry.
Consider the sequence of probability measures
$$\mu_n:={1\over \#Q_n} \sum_{a\in Q_n} \delta_a,$$
where $\delta_a$ stands for the Dirac mass at $a$. 

\begin{problem} \label{prob_eq_per}
Under suitable hypotheses on $f$, show that $\mu_n$ converges weakly to a canonical invariant probability measure $\mu$ of $f$. More precisely, for any continuous or smooth test function $\phi$ on $X$, we have
$$\lim_{n\to\infty} \langle \mu_n,\phi\rangle = \langle \mu,\phi\rangle.$$
\end{problem}

Note that if the last convergence holds and if $\mu$ has no mass on proper analytic subsets of $X$, then periodic points of $f$ are Zariski dense in $X$. More generally, the "fatness" of $\mu$ in terms of pluripotential theory is known 
 in several situations presented in the next sections.  
 
A strategy to solve the last problem is to solve Problem 1a or 1b and then to construct a good family $Q_n^{good}$ of isolated periodic points,  using tools from dynamics or complex analysis. This means $\#Q_n^{good}$ is close to an upper bound of  $\#Q_n$ and the points in $Q_n^{good}$ are equidistributed along a law canonically associated to the dynamical system.
Such a construction will  automatically give us a good lower bound for $\#Q_n$. So the question of lower bound of $\#Q_n$ is currently not a priority. Here is the type of statement that we will discuss in detail in the next sections.

\begin{theorem} \label{th_main}
Let $f$ be a dominant meromorphic correspondence on a compact K\"ahler manifold $X$. Let $Q_n$ be the set of isolated periodic points of period $n$ of $f$. Then
\begin{enumerate}
\item The cardinality of $Q_n$ grows at most exponentially fast with $n$. Moreover, its exponential rate of growth is bounded from above by the algebraic entropy of $f$, which is finite.
\item In many cases, the points in $Q_n$ are asymptotically equidistributed with respect to a canonical invariant probability measure $\mu$, i.e. the sequence $\mu_n$ converges to $\mu$ when $n$ goes to infinity.
\end{enumerate}
\end{theorem}
 
Note that even when $f$ is a holomorphic map, the first statement is new. In general, when $f$ is a meromorphic map, it is  not even  continuous and it is not obvious that the entropy is finite. On the other hand, a surprising recent result by Kaloshin says that, for smooth non-holomorphic maps, $\#Q_n$ may grow arbitrarily fast. His result contradicts the philosophy that the growth of dynamical objects should be controlled by entropy, which is finite for smooth maps. See \cite{Kaloshin} for details. However, we still believe that the philosophy holds for holomorphic, meromorphic maps or correspondences.

A main tool in complex dynamics in higher dimension is pluripotential theory, started by Lelong-Oka in 1950's. This is a powerful tool in complex analysis and geometry. The main objects in this theory are the positive closed currents which are generalizations of effective algebraic cycles in algebraic geometry. The theory is well developed for currents of bidegree $(1,1)$, i.e. currents of hypersurface type. So the first works using pluripotential theory in several  complex variables dynamics were mostly using convergence theorems for plurisubharmonic functions and their strong compactness properties.

However, unless we assume strong hypotheses on $X$ or  $f$, we need to work with positive closed currents of arbitrary bidegree. Simultaneously with complex dynamics, which is becoming quite geometric, 
the pluripotential theory is also developed. 
One of our goals is to elaborate a calculus on spaces of positive closed currents which are compact and contain the analytic cycles of arbitrary degree.
We have to consider the analogue of movable cycles and give estimates on the dimension excess  phenomenon, i.e higher dimensional 
intersection multiplicity. These questions are unavoidable even to understand the distribution of periodic points of polynomial automorphisms of $\C^k$. So we need mostly to study the geometric theory for currents which 
 is of an independent interest.  

In this paper, we will focus our discussion on Theorem \ref{th_main} and other related equidistribution problems. We will avoid technical points and emphasize the conceptual difficulties and ideas to overcome them, together with the necessary development in the theory of currents.  
Basic notions and techniques will be introduced in the next three sections. In particular, we have seen that the considered 
problems are related to the cohomology classes of the graphs of $f^n$ and hence to the action of $f^n$ on the cohomology on $X$. Crucial invariants like dynamical degrees and algebraic entropy will be introduced. 

To see the link between the distribution of orbits of varieties and the distribution of isolated periodic points, let us mention roughly a new strategy to solve Problem \ref{prob_eq_per}. The first step is to consider the current of integration on the graph $\Gamma_n$ of $f^n$, that we denote by $[\Gamma_n]$. In the interesting cases, the volume of $\Gamma_n$ grows to infinity exponentially fast. We will normalize $[\Gamma_n]$ by dividing it by a suitable constant $d_n$ and try to show that $d_n^{-1}[\Gamma_n]$ converges to a positive closed current $\T$ in $X\times X$. Imagine now that $\T$ is a generalized variety of dimension $k$. Its intersection with $\Delta$ is expected to have zero dimension. The next step is to show that there is no dimension excess for the last intersection. 
But here, we face a deep problem "how to see the dimension excess of the intersection of a current with a variety or more generally with another current"? A solution to this problem is given by our theory of densities of currents that we will briefly report on in Section \ref{section_current}. The final step is, roughly, to show that 
$$\lim_{n\to\infty} d_n^{-1} \big([\Gamma_n] \wedge [\Delta]\big) = \big( \lim_{n\to\infty} d_n^{-1} [\Gamma_n] \big)\wedge [\Delta].$$

The left hand side of the last identity is the intersection in the sense of currents and it roughly represents the probability measure $\mu_n$ equidistributed on $Q_n$. The situation will be more delicate when $\Gamma_n\cap\Delta$ is of positive dimension. We should "extract" from the intersection the part of the right dimension. The right hand side is the intersection of $\T$ with $\Delta$ and we expect to get a canonical invariant probability measure of $f$ when we identify $\Delta$ with $X$ in the canonical way. Proving the last identity is a very difficult problem in general and uses heavily the dynamical properties of the sequence $\Gamma_n$.

Now, consider the correspondence $F(x_1,x_2):=(f(x_1),f^{-1}(x_2))$ on $X\times X$. We can check that 
$$\Gamma_n=F^{-n/2}(\Delta),$$
at least in a Zariski open set. The case where $n$ is odd is just a simple technical issue. Therefore, the convergence in the first step of our approach is directly related to the equidistribution of the orbits of varieties under the action of $F^n$. Such a convergence can be proved using similar results for the correspondence $f$ on $X$. 

The equidistribution of the orbits of varieties will be presented in Section \ref{section_varieties}. 
Finally, in the last section, a list of open problems is given. Some of them require new ideas or new techniques.

\section{Meromorphic maps, correspondences and algebraic stability} \label{section_maps}

Let $(X,\omega)$ be a compact K\"ahler manifold of dimension $k$ as above. Let $\pi_1$ and $\pi_2$ denote the canonical projections from $X\times X$ to its two factors. Consider an effective $k$-cycle $\Gamma=\sum \Gamma_i$, where  $\Gamma_i$'s are irreducible analytic sets of dimension $k$ in $X\times X$. We only consider here finite sums and the $\Gamma_i$'s are not necessarily distinct. 
We assume that the restriction of $\pi_1$ to each  $\Gamma_i$ is surjective.

\begin{definition}\rm \label{def_corr}
Let $\Gamma$ be as above. We say that $\Gamma$ defines a {\it meromorphic correspondence or self-correspondence} $f$ on $X$ and $\Gamma$ is its {\it graph}. More precisely, if $A$ is any subset of $X$, we define
$$f(A):=\pi_2(\pi_1^{-1}(A)\cap\Gamma) \quad \text{and}\quad f^{-1}(A):=\pi_1(\pi_2^{-1}(A)\cap\Gamma).$$
We say that $f$ is {\it dominant} if the restriction of $\pi_2$ to each $\Gamma_i$ is surjective. In this case,
we also say that $f^{-1}$ is the {\it adjoint correspondence} of $f$. 
\end{definition}

Note that when $f$ is dominant,  $f^{-1}$ is also a meromorphic correspondence on $X$. Its graph is symmetric to $\Gamma$ with respect to the diagonal $\Delta$ of $X\times X$, i.e. the image of $\Gamma$ by the involution $(x_1,x_2)\mapsto (x_2,x_1)$. 

\begin{definition} \rm
Let $\Gamma$ and $f$ be as above.
When $\pi_1$ restricted to $\Gamma$ is generically 1:1, then $f$ is a {\it meromorphic map} and when it is 1:1, $f$ is a {\it holomorphic map}. 
\end{definition}

From now on, the correspondences we consider are all  dominant. We introduce also two indeterminacy sets
$$I(f):=\{x\in X,\ \dim f(x)>0\} \quad \text{and}\quad I(f^{-1}):=\{x\in X,\ \dim f^{-1}(x)>0\}.$$
They are respectively the {\it indeterminacy sets} for $f$ and $f^{-1}$. It is not difficult to see that $\dim \pi_1^{-1}(I(f))\cap\Gamma<k$. It follows that $\dim I(f)\leq k-2$ and a similar property holds for $I(f^{-1})$. 
If $f$ is a dominant meromorphic map such that $I(f)=\varnothing$, then $f$ is a holomorphic map. In this case, we can show that $I(f^{-1})$ is also empty but the property is true neither for correspondences nor for meromorphic maps.  The reader can find in the papers by Oguiso and Truong  \cite{Oguiso1, OT} some recent examples of dynamically interesting meromorphic maps.

Consider now two dominant meromorphic correspondences $f$ and $f'$ on $X$. We can define the correspondence $f'\circ f$ in the following way. Choose a Zariski open set $\Omega$ in $X$ such that for $x\in\Omega$ the set $f(x)$ is finite. It is not difficult to see that $f(\Omega)$ is a Zariski open set of $X$. By reducing $\Omega$, we can also assume that $f'(x)$ is finite for $x\in f(\Omega)$. 
We will define $(f'\circ f)(x)$ as the set  $f'(f(x))$ for $x\in \Omega$. The graph of this multi-valued map is an effective cycle in $\Omega\times X$. 
Let $\widetilde\Gamma$ be the closure of this cycle in $X\times X$ which is, by definition,  the graph of the dominant correspondence $f'\circ f$ on $X$. 

The construction of $f'\circ f$ can be obtained in a more geometrical way. Let $\Gamma$ and $\Gamma'$ denote the graphs of $f$ and $f'$ in $X\times X$. Consider the product $\Gamma\times\Gamma'$ in $X^4$ and let $(x_1,x_2,x_3,x_4)$ denote a point in $X^4$. Consider the projection $\widehat\Gamma$ of the intersection $(\Gamma\times\Gamma')\cap \{x_2=x_3\}$ into the product of the first and the last factors in $X^4$.  We obtain $\widetilde \Gamma$  from $\widehat\Gamma$ by removing components of dimension larger than $k$ and components whose projections onto the factors of $X\times X$ are not surjective.
This is the graph of the composition $f'\circ f$, see Definition \ref{def_corr}. 

We can define the iterate $f^n:= f\circ \cdots\circ f$, $n$ times, for every $n\geq 1$. 
Denote also by $f^{-n}$ the adjoint of $f^n$ which is also the iterate of order $n$ of the correspondence $f^{-1}$. We will discuss now the notion of algebraic stability of $f$ which plays an important role in our study.

Let $S$ be a $(p,q)$-current on $X$. We define formally the {\it pull-back} of $S$ by $f$ by
$$f^*(S):=(\pi_1)_* (\pi_2^*(S)\wedge [\Gamma]),$$
when the last expression makes sense. Note that the operators $\pi_i^*$ and $(\pi_i)_*$ are well-defined on all currents. Therefore, the last definition is meaningful when the wedge-product $\pi_2^*(S)\wedge [\Gamma]$ is meaningful. We also define the {\it push-forward} operator $f_*$ as the pull-back operator $(f^{-1})^*$ associated to $f^{-1}$.

Consider the particular case of a continuous or smooth differential $(p,q)$-form $\phi$ on $X$. The wedge-product $\pi_2^*(\phi)\wedge [\Gamma]$ is well-defined because $\pi_2^*(\phi)$ is a continuous form and $[\Gamma]$ is a current of order 0. 
So $f^*(\phi)$ is well-defined in the sense of currents. However, the value of $f^*(\phi)$ at a point $x$ is roughly the sum of the values of $\pi_2^*(\phi)$ on the fiber $\pi_1^{-1}(x)\cap\Gamma$. We can check that $f^*(\phi)$ is in general an $L^1$ form but it may be singular at the critical values of the map $\pi_1$ restricted to $\Gamma$ which contain the indeterminacy set $I(f)$. So we cannot iterate the operator $f^*$ on continuous or smooth forms : for example, the expression $f^*(f^*(\phi))$ is not meaningful in general. So it is necessary to establish a calculus.

Recall that the Hodge cohomology group $H^{p,q}(X,\C)$ of $X$ can be defined using either smooth forms or singular currents. Observe that when $\phi$ is closed or exact then $f^*(\phi)$ is also closed or exact. Therefore, the above operator $f^*$ induces a linear map from $H^{p,q}(X,\C)$ to itself, that we still denote by $f^*$.  The operator $f_*$ on $H^{p,q}(X,\C)$ is defined to be $(f^{-1})^*$. We can iterate those operators as for every linear operator on a vector space. The following definition is a  slight extension of a notion introduced by Forn\ae ss and the second author, see \cite{FS1}.

\begin{definition} \rm
We say that a meromorphic correspondence  $f$ is {\it algebraically $p$-stable}, $0\leq p\leq k$, if we have $(f^n)^*=(f^*)^n$ on $H^{p,p}(X,\C)$ for every $n\geq 1$. 
We say that $f$ is {\it algebraically stable} if it is algebraically $1$-stable and $f$ is {\it totally algebraically stable} if it is algebraically $p$-stable for every $p$.
\end{definition}

Note that $f$ is always algebraically $p$-stable for $p=0$ and $p=k$. If $f$ is a dominant holomorphic map then it is totally algebraically stable. There exist meromorphic maps which are not algebraically stable. 
For example, in dimension $k=2$, if a curve is sent to an indeterminacy point then the map is not algebraically stable.
In some situations, the algebraic stability can be easily checked but in general this seems to be a difficult question since we need to check the identity $(f^n)^*=(f^*)^n$ for all $n$, see also Nguyen \cite{Nguyen}.

Recall that in the construction of the graph of the composition of two correspondences, one step is to eliminate bad components of some analytic cycles. This step is in fact the cause of the algebraic instability. 
The reader can observe the phenomenon with the rational involution $f$ on $\P^2$ given on $\C^2$ by  $f(z_1,z_2)=(1/z_1,1/z_2)$.
The algebraic instability makes the dynamical study of $f$ much more difficult. Therefore, in many results, the algebraic stability is assumed. It is  verified for large classes of maps. For example, in a projective space, the set of algebraically unstable maps of a fixed degree form a countable union of analytic sets in the parameter space. 

Finally, observe that if $\varphi$ and $\psi$ are smooth differential forms of bidegrees $(p,q)$ and $(k-p,k-q)$ respectively, we have 
$$\langle f^*(\phi), \psi\rangle =\langle \phi, f_*(\psi)\rangle.$$
It follows that the operator $f^*:H^{p,q}(X,\C)\to H^{p,q}(X,\C)$ is dual to the operator $f_*:H^{k-p,k-q}(X,\C)\to H^{k-p,k-q}(X,\C)$ via the Poincar\'e's duality. So the algebraic stability can be also expressed in terms of $(f^n)_*$.

\section{Positive closed currents, super-potentials and densities} \label{section_current}

In this section, we will report on some recent techniques used to deal with positive closed currents of arbitrary bidegree. We refer the reader to \cite{Demailly, Hormander, Siu} for basic notions and results of pluripotential theory and to \cite{Voisin} for Hodge theory. 
Positive closed currents can be seen as positive closed differential forms with distribution coefficients. This is a common point of view in the analytic aspect. Locally, positive closed currents can be approximated by smooth positive closed forms using the standard process of convolution. Computation with positive closed currents is possible thanks to an appropriate control of the regularization  process. 

On a compact K\"ahler manifold, it is more convenient to work in the global setting rather than in the local setting. There are two reasons. The first one is that one often use the integration by parts or more generally the Stokes formula. Working in the global setting allows us to avoid boundary terms. The second reason is even more important. In a compact K\"ahler manifold $(X,\omega)$ of dimension $k$, if $T$ is a positive closed $(p,p)$-current, the pairing $\langle T,\omega^{k-p}\rangle$ depends only on the (Hodge or de Rham) cohomology classes of $T$ and of $\omega$. This quantity is equivalent to the mass of $T$ which is, by definition, the norm of $T$ as a linear operator on the space of continuous $(k-p,k-p)$-forms. Therefore, a large part of computations with positive closed currents and even for positive $\ddc$-closed currents reduces to a computation with cohomology classes which is often simpler.
 
Except for homogeneous manifolds, the convolution process doesn't work in the global setting. Moreover, in general, we cannot regularize positive closed currents without loosing the positivity.
The following result gives us a regularization with a control of the positivity loss, see \cite{Demailly, DS1}. For convenience, we also call   
$\|T\|:=\langle T,\omega^{k-p}\rangle$ {\it the mass} of $T$. 

\begin{theorem}[Demailly for $p=1$, Dinh-Sibony for $p\geq 1$] \label{th_reg}
Let $(X,\omega)$ be a compact K\"ahler manifold. There is a constant $c>0$ depending only on $X$ and $\omega$ satisfying the following property. If $T$ is a positive closed $(p,p)$-current on $X$, there are positive closed $(p,p)$-currents $T^+$ and $T^-$ which can be approximated by smooth positive closed $(p,p)$-forms and such that 
$$T=T^+-T^- \qquad \text{and} \qquad \|T^\pm\|\leq c \|T\|.$$   
\end{theorem}

The theorem also holds for other classes of currents, e.g. positive $\ddc$-closed currents. It is similar to the known fact in algebraic geometry that any cycle can be written as the difference of movable effective cycles.
The regularization process hidden in the last theorem preserves good properties of $T$ when they exist. For example, if $T$ is smooth in some open set $U$, then $T^\pm$ are also smooth there and the approximation of $T^\pm$ by smooth positive closed forms on $X$ is uniform on compact subsets of $U$. Specific needs can be obtained by going through the details of the proof of the above theorem, see \cite{Demailly, DNT1, DS1}.

We will discuss now the notion of super-potentials.  The starting point is that the pluripotential theory is well developed for positive closed currents of bidegree $(1,1)$ thanks to the notion of plurisubharmonic (p.s.h. for short) functions. More precisely, if $T$ is a positive closed $(1,1)$-current, then we can write locally $T=\ddc u$ in the sense of currents, where $u$ is a p.s.h. function and $\dc:={i\over 2\pi} (\partial -\dbar)$. Working with a pointwise defined function is more confortable than working directly with a current. 
For example,  it allows more operations like multiplying a positive current $S$ by $u$. One just has to assume that $u$ is integrable with respect to the trace measure of $S$. But when $T$ is of higher bidegree, the potentials is just an $L^1$-form, one cannot consider their wedge-product with $S$ if $S$ is singular.

In the global setting, if $\alpha$ is a real smooth closed $(1,1)$-form in the cohomology class of $T$, we can write globally 
$$T=\alpha+\ddc u,$$ 
where $u$ is a quasi-p.s.h. function, i.e., locally the sum of a p.s.h. function and a smooth one. This function $u$ is unique up to an additive constant. So if we normalize $u$ by the condition 
$$\int_X u \omega^k=0$$
then $u$ is unique. We call $u$ the {\it normalized quasi-potential} of $T$. One of the key technical point in the use of quasi-p.s.h. functions is the control of their singularities. The 
most useful general property is perhaps the following consequence of a theorem by Skoda.

\begin{theorem}[Skoda] \label{th_Skoda}
Let $X,\omega$ and $\alpha$ be as above. Then there are constants $c>0$ and $\lambda>0$ such that if $T$ is a positive closed $(1,1)$-current in the cohomology class of $\alpha$ and $u$ is its normalized quasi-potential, then we have
$$\int_X e^{\lambda|u|} \omega^k\leq c.$$
\end{theorem} 

It is easy to deduce estimates of the  $L^p$-norm of $u$ for all $1\leq p<\infty$. Refined versions of this theorem can be found in \cite{DNS1, Kaufmann, Vu}. So quasi-p.s.h. functions are almost as good as bounded functions. When the above estimate is verified for a probability measure $\nu$, in place of $\omega^k$, we say that the measure $\nu$ is moderate. The support of a moderate measure is Zariski dense.

Super-potentials are functions which play the role of quasi-potentials for positive closed currents of arbitrary bidegree. We will introduce them briefly and refer the reader to \cite{DS3, DS4} for details. Let $\alpha$ be a real smooth closed $(p,p)$-form. Let $T$ be a positive closed $(p,p)$-current in the cohomology class of $\alpha$. Let $\Cc_{k-p+1}(X)$ denote the set of positive closed currents of bidegree $(k-p+1,k-p+1)$ and mass 1 in $X$ and $\widetilde\Cc_{k-p+1}(X)$ the set of currents in $\Cc_{k-p+1}(X)$ which are smooth forms. The normalized super-potential of $T$ is a real-valued function which is defined at least on $\widetilde\Cc_{k-p+1}(X)$ and may be extended to a function in a larger set of currents. It is denoted by $\Uc_T$ and given by 
$$\Uc_T(R):=\langle U_T,R\rangle \quad \text{for} \quad R\in \widetilde\Cc_{k-p+1}(X),$$
where $U_T$ is a normalized quasi-potential of $T$, that is, $U_T$ is a $(p-1,p-1)$-current such that 
$$\ddc U_T=T-\alpha \quad \text{and} \quad \langle U_T,\omega^{k-p+1}\rangle=0.$$

Note that when $p>1$ the quasi-potential $U_T$ also exists but is not unique. However, we can show that the definition of $\Uc_T$ does not depend on the choice of $U_T$. So the super-potential $\Uc_T$ is a canonical function associated to $T$ and $\alpha$. It enjoys several properties similar to quasi-p.s.h. functions but we don't quote all of them here. Since it is defined on a space of infinite dimension, it is not clear how to get a similar property as the above Skoda's estimate for quasi-p.s.h. functions. The following result gives the answer to this question, see \cite{DS3,DS4}.

\begin{theorem}[Dinh-Sibony] \label{th_super_pot}
Let $X,\omega$ and $\alpha$ be as above. There is a constant $c>0$ such that if $T$ is a positive closed $(p,p)$-current in the cohomology class of $\alpha$ and $\Uc_T$ is its normalized super-potential, then 
$$|\Uc_T(R)|\leq c(1+\log^+\|R\|_{\Cc^1}),$$
where $\log^+:=\max(\log,0)$.
\end{theorem} 

Super-potentials also allow to define the intersection of positive closed currents of arbitrary bidegree. For example, if the super-potential $\Uc_T$ of $T$ can be extended to a continuous function on $\Cc_{k-p+1}(X)$, then for any positive closed $(q,q)$-current $S$ in $X$ with $1\leq q\leq k-p$, using the theory one shows that the wedge-product $T\wedge S$ is well-defined and depends continuously on $S$. This notion was used in dynamics to define invariant measures as the intersection of positive closed invariant currents. We refer to \cite{DS3,DS4,Vu2} for details.

In the rest of this section, we will briefly introduce the theory of densities for positive closed currents, see \cite{DS5} for details. As we have seen in the introduction, the study of periodic points of a map or correspondence is strongly related to the so called dimension excess of the intersection of analytic cycles. If $V_1$ and $V_2$ are two analytic subsets of $X$ of co-dimension respectively $p_1$ and $p_2$, then their intersection is expected to be empty when $p_1+p_2>k$. Otherwise, if their intersection is non-empty, its dimension is expected to be $k-p_1-p_2$. In general, this dimension is at least equal to $k-p_1-p_2$. However, in the first case, the intersection may be non-empty and in the second case the dimension may be larger than $k-p_1-p_2$. We refer to this well-known phenomenon as the dimension excess. 

The starting point of the theory of densities for positive closed currents is to determine if such a phenomenon happens for general positive closed currents. The theory is crucial in dynamics because even when we want to understand algebraic cycles we cannot limit ourself to cycles of bounded degree  and therefore positive closed currents appear naturally as limits of normalized cycles.
It is also convenient to work with the space of positive closed currents of bidegree $(p,p)$ because it admits good compactness properties. 

There exists a particular case for which we first explain our point of view.  This is the notion of Lelong number. If $T$ is a positive closed $(p,p)$-current in $X$ and $a$ is a point in $X$, then the Lelong number $\nu(T,a)$ of $T$ at $a$ represents the density of $T$ at the point $a$. In our point of view, when this number is strictly positive, the "intersection" of $T$ with the point $a$ has a dimension excess. When $T$ is given by a subvariety of $X$, this phenomenon appears only when $a$ belongs to this variety. Lelong number was originally defined using local holomorphic coordinates near $a$. The following result is fundamental in the theory.

\begin{theorem}[Lelong, Siu]
Let $X,T$ and $a$ be as above. Then 
\begin{enumerate}
\item The Lelong number $\nu(T,a)$ is intrinsic, i.e. it does not depend on the choice of local coordinates.
\item The Lelong number $\nu(T,a)$ is upper semi-continuous with respect to $T$. In particular, it is bounded by a constant times the mass of $T$.
\item For every $c>0$, the upper level set $\{a\in X, \nu(T,a)\geq c\}$ is an analytic subset of $X$.
\end{enumerate} 
\end{theorem}

We describe now the case where the point $a$ is replaced by a submanifold $V$ of dimension $l$ of $X$. The aim is to see if there is a dimension excess for the "intersection" between $T$ and $V$. Let $N_{V|X}$ be the normal vector bundle of $V$ in $X$. We identify its zero section with the manifold $V$.  Let $\tau$ denote a smooth map from a neighbourhood of $V$ in $X$ to a neighbourhood of $V$ in $N_{V|X}$ such that the restriction of $\tau$ to $V$ is the identity and the induced map from $N_{V|X}$ to itself is also identity. Denote by $\overline N_{V|X}$ the natural compactification of $N_{V|X}$ and by $A_\lambda : N_{V|X}\to N_{V|X}$, for $\lambda\in \C^*$, the map which is the multiplication by $\lambda$ in each fiber of $N_{V|X}$. This map extends to an automorphism of $\overline N_{V|X}$.

Let $T$ be a positive closed $(p,p)$-current. For simplicity, assume that it has no mass on $V$. 
Consider now the family of currents $T_\lambda:=(A_\lambda)_*\tau_*(T)$ in open sets of $\overline N_{V|X}$. 
The domain of definition of $T_\lambda$ increases to $N_{V|X}$ when $|\lambda|$ tends to infinity.
In general, these currents are not 
positive nor of bidegree $(p,p)$. However, we obtain the following result. 

\begin{theorem}[Dinh-Sibony] \label{th_density}
Let $X,V, \overline N_{V|X}, \tau, A_\lambda$ and $T$ be as above.
\begin{enumerate}
\item The family $\{T_\lambda\}$ is relatively compact in the sense that for any sequence $\lambda_n\to\infty$, there is a subsequence $\lambda_{n_i}$ such that $T_{\lambda_{n_i}}$ converges to a current $T_\infty$ on $N_{V|X}$.
\item The current $T_\infty$ does not depend on the choice of the map $\tau$. Moreover, it is a positive closed $(p,p)$-current which can be extended by $0$ to a positive closed $(p,p)$-current on $\overline N_{V|X}$, that we still denote by $T_\infty$. The mass of the last current is bounded by a constant times the mass of $T$.
\item The cohomology class  of $T_\infty$ in $H^{p,p}(\overline N_{V|X},\C)$, denoted by $\kappa^V(T)$, is intrinsic, i.e. it does not depend on the choice of $\tau$ nor on the choice of the sequence $\lambda_n$. Moreover, it is upper semi-continuous with respect to $T$. 
\end{enumerate}
\end{theorem}

The class $\kappa^V(T)$ is called the {\it total tangent class} of $T$ along $V$. It represents the density of $T$ near $V$.
Note that the bi-dimension of $T$ is $(k-p,k-p)$. So if $p>l$, we expect, at least in the generic case, that this class vanishes or equivalently $T_\infty=0$. In other words, there is no dimension excess for the "intersection" between $T$ and $V$. When this class is not zero, we conclude that there is a dimension excess. This situation corresponds to the case of positive Lelong number when $V$ is reduced to a point. 

Consider now the case where $p\leq l$. It is a little bit more subtle  to observe the phenomenon of dimension excess. 
By K\"unneth theorem, the Hodge cohomology class  $\kappa^V(T)$ can be identified to a polynomial 
with coefficients $\kappa_j^V(T)$ in 
$H^{l-j,l-j}(V,\C)$ for $l-p\leq j\leq l$.
Let $s$ denote the maximal number such that $\kappa_s^V(T)\not=0$. We take $s=l-p$ if all $\kappa_j^V(T)$ vanish.  This number is called the {\it tangential h-dimension} of $T$ along $V$. In the generic situation, we expect it to be minimal, i.e. equal to $l-p$. 
When this dimension is larger than $l-p$, we conclude that there is a dimension excess for the intersection between $T$ and $V$. 

We can define the densities between two positive closed currents $T$ and $S$ by considering the densities of $T\otimes S$ along the diagonal $\Delta$ in $X\times X$. Note also that
the theory of densities of currents also allows to define intersection of positive closed currents of arbitrary bidegree under suitable conditions, including the absence of the dimension excess. This is a promising research direction which may offer the most general setting where we can define the wedge-product of currents. The reader will find in \cite{DS5} a more detailed exposition.

\section{Dynamical degrees, topological and algebraic entropies} \label{section_degree}

In complex dynamics, a crucial question is to understand the action of the map or correspondence on positive closed currents. This action, if it is well-defined, is often compatible with the action on cohomology and permits to deal with the distribution of orbits and to construct dynamically significant invariant measures.
Moreover, the mass of a positive closed current depends only on its cohomology class. 
Therefore, in order to control the mass of these currents, 
it is very useful to understand the action of the map or correspondence on cohomology. 
We will give in this section some general basic properties.
Note that although we work with linear actions on  finite dimensional vector spaces, it is a difficult problem to describe their behaviour when we iterate the map or correspondence because of the algebraic instability mentioned in Section \ref{section_maps}. 

Let $(X,\omega)$ be a compact K\"ahler manifold of dimension $k$ as above. Let $f$ be a dominant meromorphic correspondence on $X$. We also fix some norms on the Hodge cohomology groups.

\begin{definition} \rm
We call {\it dynamical degree} of order $p$ of $f$ the following limit
$$d_p(f):=\lim_{n\to\infty} \|(f^n)^*: H^{p,p}(X,\C)\to H^{p,p}(X,\C)\|^{1/n}$$
and {\it algebraic entropy} of $f$ the following constant
$$h_a(f):=\max_{0\leq p\leq k} \log d_p(f).$$ 
The last topological degree $d_k(f)$ is also called {\it topological degree} because it is equal to the number of points in $f^{-1}(a)$ for a generic point $a$ in $X$.
\end{definition}

Note that the discussion at the end of Section \ref{section_maps} implies that 
$$d_p(f):=\lim_{n\to\infty} \|(f^n)_*: H^{k-p,k-p}(X,\C)\to H^{k-p,k-p}(X,\C)\|^{1/n}=d_{k-p}(f^{-1}).$$
Therefore, we also have 
$$h_a(f)=h_a(f^{-1}).$$
We have the following general result.

\begin{theorem}[Dinh-Sibony]
The limit in the above 
definition of $d_p(f)$ always exists. It is finite and doesn't  depend on the choice of the norm on $H^{p,p}(X,\C)$. 
Moreover, the dynamical degrees and the algebraic entropy are bi-meromorphic invariants of the dynamical system. That is, if $\pi:X'\to X$ is a bi-meromorphic map from a compact K\"ahler manifold $X'$ to $X$, then 
$$d_p(\pi^{-1}\circ f\circ \pi)=d_p(f) \quad  \text{and}  \quad h_a(\pi^{-1}\circ f\circ \pi)=h_a(f).$$
We also have for $n\geq 1$ that
$$d_p(f^n)=d_p(f)^n  \quad  \text{and}  \quad h_a(f^n)=n h_a(f).$$
\end{theorem}

Note that when $X$ is a projective space, the first statement was used by Forn\ae ss and the second author for $p=1$ in order to construct the Green dynamical $(1,1)$-current \cite{FS1}. Also for projective spaces, it was extended by Russakovskii-Shiffman for higher degrees \cite{RS}. In this case, the group $H^{p,p}(X,\C)$ is of dimension 1 and the action of $(f^n)^*$ is just the multiplication by a constant $d_{p,n}$. Therefore, we easily get $d_{p,n+m}\leq d_{p,n}d_{p,m}$ which implies the result.
The proof of the above theorem in the general case uses in an essential way a computation with positive closed currents and Theorem \ref{th_reg} plays a crucial role. We refer to \cite{DNT1, DS1, DS12, Truong1, Truong2} for details and  some extensions of this result. 
We also obtained in these works the following result, which has been obtained by Gromov for holomorphic maps \cite{Gromov2}.

\begin{theorem}[Gromov, Dinh-Sibony]
Let $X$ and $f$ be as above. Then  the topological entropy $h_t(f)$ of $f$ is bounded from above by its algebraic entropy $h_a(f)$.
In particular, the topological entropy of $f$ is finite.
\end{theorem}

The topological entropy measures the rate of divergence of the orbits of points. Its definition for meromorphic maps or correspondences is the same as the Bowen's definition for continuous maps, except that we don't consider orbits which reach the indeterminacy set. Therefore, it is not obvious that the entropy of a meromorphic map is finite.  
Note also that when $f$ is a holomorphic map, the above result  combined with a theorem by Yomdin \cite{Yomdin} implies that the topological entropy is indeed equal to the algebraic one. We expect that such a property still holds for large families of meromorphic maps and correspondences. We don't know if in general, there is always a map or correspondence $\hat f$ bi-meromorphically conjugate to $f$ such that $h_t(\hat f)=h_a(\hat f)$. If $X$ is a projective manifold, by composing holomorphic projections from $X$ to $\P^k$ with the adjoints of such maps, it is easy to construct correspondances of positive topological entropy on $X$.

Observe that the action of $f^n$ on $H^{p,q}(X,\C)$ is not explicitly used in the above property of entropies when $p\not=q$. This can be explained by the following inequality obtained by the first author in \cite{Dinh1}
$$\limsup_{n\to\infty} \|(f^n)^*: H^{p,q}(X,\C)\to H^{p,q}(X,\C)\|^{1/n}\leq \sqrt{d_p(f)d_q(f)}.$$

Dynamical degrees are not easy to compute except in the case of algebraically $p$-stable correspondences. The following result was obtained by the authors for automorphisms  in \cite{DS13}  and for meromorphic maps in  unpublished lecture notes. It implies that for a holomorphic family of algebraically $p$-stable correspondences the dynamical degree of order $p$ is constant.

\begin{proposition}
Let $f$ be an algebraically $p$-stable correspondence on a projective manifold $X$. 
Then the dynamical degree $d_p(f)$ of $f$ takes values in a discrete subset of algebraic integers in $[0, +\infty)$.
\end{proposition}

Let $\NS_p(X,\R)$ be the N\'eron-Severi subspace of $H^{p,p}(X,\C)$ spanned by the classes of algebraic $(k-p)$-cycles in $X$. Then the action of $f$ on $\NS_p(X,\R)$ preserves the lattice spanned by the algebraic $(k-p)$-cycles with integer coefficients. 
Since $f$ is algebraically stable, we can show that $d_p(f)$ is the spectral radius of the action $f^*$ on $\NS_p(X,\R)$ and therefore is the largest root of a monic polynomial with integer coefficients. The degree of this polynomial is the dimension of $\NS_p(X,\R)$. The proposition follows easily. When the map or correspondence is not algebraically $p$-stable, the following problem seems to be difficult.

\begin{problem}
Let $f$ be an arbitrary dominant meromorphic map or correspondence on a compact K\"ahler manifold. Are its dynamical degrees always algebraic integers ?
\end{problem}

We continue our discussion on dynamical degrees.
Recall that a direct consequence of the mixed Hodge-Riemann theorem applied to the graphs of $f^n$, see e.g. \cite{DN1, Gromov1}, implies that, when $f$ is a meromorphic map,  the function $p\mapsto \log d_p(f)$ is concave. Equivalently, we have
$$ d_p(f)^2\geq d_{p-1}(f)d_{p+1}(f) \quad \text{for} \quad 1\leq p\leq k-1.$$
In particular, we have $1\leq d_p(f)\leq d_1(f)^p$, $h_a(f)>0$ if and only if $d_1(f)>1$, and there are two numbers $r$ and $s$ with $0\leq r\leq s\leq k$ such that
$$1=d_0(f)<\cdots<d_r(f)=\cdots=d_s(f)>\cdots>d_k(f).$$
The property still holds for a correspondence $f$ when the graphs of  $f^n$ are irreducible.

If $V$ is an analytic subset of pure dimension $k-p$ in $X$, the compactification of $f(V\setminus I(f))$ is also an analytic subset of pure dimension $k-p$. We call it the strict transform of $V$ by $f$. Similarly, if $T$ is a positive closed $(p,p)$-current on $X$, we can define the push-forward of $T$ by $f$ restricted to $X\setminus I(f)$. Using Theorem \ref{th_reg}, one can show that the obtained current can be extended by 0 to a positive closed $(p,p)$-current on $X$, see \cite{DS2} for details. This is  the so-called strict transform of $T$ by $f$.

Let $V'$ be an analytic subset of pure dimension $p$ and $T'$ a positive closed $(k-p,k-p)$-current on $X$. Let $\epsilon>0$ be any constant. It is possible to prove, using again Theorem \ref{th_reg}, that the $2p$-dimensional volume (resp. the mass) of the strict transform of $V'$ (resp. $T'$) by $f^{n}$ is smaller than a constant times $(d_p(f)+\epsilon)^n$. Moreover, $d_p(f)$ is the smallest constant satisfying this property. Therefore, if $r$ and $s$ are as above, the dynamical system  $f$ has roughly $r$ expanding directions, $s-r$ neutral directions and $k-s$ contracting directions.
In dynamics, neutral directions often make the study more difficult.

\begin{definition} \rm
Let $f$ be a dominant meromorphic correspondence on $X$. We say that $f$ is {\it algebraically hyperbolic} if it admits a dynamical degree $d_p(f)$ which is strictly larger than the other ones. If moreover, this is the last dynamical degree, we say that $f$ is {\it algebraically expanding} or {\it with dominant topological degree}. Let $f$ be algebraically hyperbolic as above and assume that
the action  of $f^*$ on $H^{p,p}(X,\C)$ has only one eigenvalue of maximal modulus which is simple and equal to $d_p(f)$.
We say that the action of $f$ on cohomology is {\it simple} or $f$ is {\it algebraically simple}, see \cite{DS4}.
\end{definition}

Note that in some references, one uses the word "cohomologically" instead of "algebraically". However, we think the second choice is more appropriate because the action of $f$ on cohomology is not hyperbolic in the usual sense because it often has eigenvalues of modulus 1. Most of results in complex dynamics are obtained under the condition that the action of the map on cohomology is simple. Note that when $f$ is algebraically expanding, then its action on cohomology is simple because $\dim H^{k,k}(X,\C)=1$.

Finally, though analytic tools are dominant in the study of the objects introduced in this section, many questions are of algebraic nature and can be asked for maps or correspondences which are defined over fields different from $\C$. The reader will find in 
\cite{EOY,ES,Truong3} and the references therein some recent algebraic counterparts of the topics presented in this paper.

\section{Number of isolated periodic points and equidistribution} \label{section_periodic}

In this section, we will discuss in detail Theorem \ref{th_main} stated in the introduction. 
We first have the following general result recently obtained in \cite{DNT3}.

\begin{theorem}[Dinh-Nguyen-Truong] \label{th_upper_bound_per}
Let $f$ be a dominant meromorphic correspondence on a compact K\"ahler manifold $X$ and let $h_a(f)$ be the algebraic entropy of $f$. If $Q_n$ denotes the set of isolated periodic points of period $n$ of $f$ counted with multiplicity, then we have
$$\limsup_{n\to\infty}{1\over n} \log \#Q_n\leq h_a(f).$$
In particular, $f$ is  an Artin-Mazur correspondence, that is, the cardinality of $Q_n$ grows at most exponentially fast with $n$.
\end{theorem}

The proof of this theorem uses the theory of densities for positive closed currents, described in Section \ref{section_current}. 
Let $[\Gamma_n]$ denote the current of integration on the graph $\Gamma_n$ of $f^n$. We can show that its mass (or equivalently the $2k$-dimensional volume of $\Gamma_n$) is $O(e^{n\lambda})$ for every constant $\lambda>h_a(f)$. It follows that the family of currents $e^{-n\lambda}[\Gamma_n]$ is bounded. So their densities along the diagonal $\Delta$ are also bounded. On the other hand, by definition, the total tangent class of $[\Gamma_n]$ along $\Delta$ is at least equal to $\#Q_n$ times the cohomology class of a fiber of the bundle $\overline N_{\Delta|X^2}$. We then deduce that $\#Q_n = O(e^{n\lambda})$, which implies the theorem.

The following result is an immediate consequence of the last theorem.

\begin{corollary}
Let $f, X$ and $Q_n$ be as in Theorem \ref{th_upper_bound_per}. Then
the zeta-function associated to $f$
$$\zeta_f(z):=\sum_{n\geq 1} {1\over n} (\#Q_n) z^n$$ 
is analytic in a neighborhood of $0\in\C$.
\end{corollary}

The following question is related to the existence of some recurrence property of the sequence $\#Q_n$. A positive answer is known in some situations where $\#Q_n$ can be computed explicitly. 

\begin{problem}
Does the zeta-function $\zeta_f$ always admit a meromorphic or rational extension to the whole complex plane $\C$ ?
We can ask the same question for the zeta-function defined with the sequence of Lefschetz numbers $L(f^n)$.
\end{problem}

We continue the discussion on Problems 1a and 1b. The theory of densities of currents  can be used to get more precise estimates. We will give here two examples. The following result was obtained in \cite{DNT2}. Recall the Lefchetz number $L(f^n)$ is the intersection number associated to $\Gamma_n$ and $\Delta$. 

\begin{theorem}[Dinh-Nguyen-Truong] \label{th_d_k_per}
Let $f$ be a meromorphic correspondence on $X$. Assume that $f$ is algebraically expanding, i.e. the last dynamical degree $d_k(f)$ of $f$ is strictly larger than the other dynamical degrees. Then we have
$$L(f^n)=d_k(f)^n+o(d_k(f))^n \quad \text{and} \quad \#Q_n \leq d_k(f)^n+o(d_k(f))^n= L(f^n) +o(L(f^n)).$$
\end{theorem}

The result was stated in \cite{DNT2} for maps but the proof is the same for correspondences. This is the solution to Problem 1b in the considered setting. We will see later that the last inequality is in fact an equality. The first identity is not difficult to obtain using the action of $(f^n)^*$ on cohomology described in Section \ref{section_degree}, in particular, the action on $H^{p,q}(X,\C)$ with $p\not=q$. 

The key point in the proof of the inequality in the theorem is that the sequence of currents $d_k(f)^{-n}[\Gamma_n]$ converges to a positive closed $(k,k)$-current $\T$ in $X\times X$ which is equal to $\pi_1^*(\mu)$ for some probability measure $\mu$ on $X$. The current $\T$ is "transverse" to the diagonal $\Delta$ and therefore there is no dimension excess between $\T$ and $\Delta$. We easily see that the total tangent class of $\T$ along $\Delta$ is equal to  the cohomology class of a fiber of the bundle $\overline N_{\Delta|X^2}$. Then, the inequality in the last theorem follows from the upper semi-continuity property in Theorem \ref{th_density}.

In the following situation, we also have a satisfactory answer to Problem 1b, see \cite{DNT2}. 

\begin{theorem}[Dinh-Nguyen-Truong]
Let $f$ be a dominant meromorphic map on a compact K\"ahler surface $X$. Assume that $f$ is algebraically stable and that its first dynamical degree $d_1(f)$ is strictly larger than its topological degree $d_2(f)$. Then we have 
$$L(f^n)=d_1(f)^n+o(d_1(f))^n \quad \text{and} \quad \#Q_n \leq d_1(f)^n+o(d_1(f))^n= L(f^n) +o(L(f^n)).$$
\end{theorem}

The proof of this theorem follows the same strategy. We show that the limit $\T$ in this case has the form $T_+\otimes T_-$ where $T_\pm$ are positive closed $(1,1)$-currents on $X$, which cannot both have mass on the same analytic curve. The last property allows us to show that there is no dimension excess for the intersection of $\T$ with $\Delta$. This is the key point in the proof. The result can be extended to meromorphic correspondences under the condition that its action on cohomology is simple.

The reader will find some related results in Favre \cite{Favre}, Iwasaki-Uehara \cite{IwasakiUehara}, Saito \cite{Saito}, Xie \cite{Xie}  and the references therein. Note that in some references, periodic points which are indeterminacy points may not be counted. 
Here we count them and they can be the only intersection points of the graphs $\Gamma_n$ and $\Delta$. This this the case for the rational extension to $\P^2$ of the polynomial map
$f(z,w)=(z+1, z^2+w)$. 

We will now consider some situations where the equidistribution of periodic points is known. The following result was obtained in the case of holomorphic endomorphisms of $\P^k$ by Briend-Duval \cite{BD1}. The general version used Theorem \ref{th_d_k_per} above was obtained by Nguyen, Truong and the first author in \cite{DNT2} for maps but the proof can be extended to correspondences.

\begin{theorem}[Briend-Duval, Dinh-Nguyen-Truong] \label{th_eq_dominant}
Let $X, f, d_k(f)$ and $Q_n$ be as in Theorem \ref{th_d_k_per}. Then there is an invariant  probability measure $\mu$ such that 
$$\lim_{n\to\infty} {1\over d_k(f)^n} \sum_{a\in Q_n} \delta_a =\mu,$$
where $\delta_a$ denotes the Dirac mass at $a$.
\end{theorem}

This result shows that isolated periodic points of order $n$ of $f$ are equidistributed with respect to $\mu$ when $n$ goes to infinity. 
It solves Problem 2 in this situation. The statement still holds if we replace $Q_n$ by the subset of repelling periodic points of period $n$. 

There are different ways to construct the measure $\mu$ : the theorem gives one way to obtain it and we can also get it as the probability measure such that the sequence of currents $d_k(f)^{-n}[\Gamma_n]$, considered above, converges to $\pi_1^*(\mu)$. 
So $\mu$ can be seen as the intersection of the limit of  $d_k(f)^{-n}[\Gamma_n]$ with the diagonal $\Delta$ of $X\times X$. 
We will see in the next section another way to produce this measure. It is known that this measure has no mass on proper analytic subsets of $X$. Therefore, periodic points are Zariski dense in $X$. Similar properties hold for other situations that we consider below.

In order to prove this theorem, using Theorem \ref{th_d_k_per}, we only need to construct a good family of periodic points. At least in the case where $X$ is a projective manifold, the strategy described in the introduction works quite well. The construction used in the original proof of this theorem follows the approach developed in \cite{DS0}, which is based on the construction of the inverse branches of $f^n$ defined on balls in $X$. The construction is quite delicate because of the needed size control. Note that the measure $\mu$ here doesn't have mass on proper analytic subsets of $X$ and satisfies $f^*(\mu)=d_k(f)\mu$.

\smallskip

Here is another simple situation where we can solve Problem 2.

\begin{theorem} \label{th_per_surface}
Let $f$ be a holomorphic automorphism of a projective surface $X$. Assume that its first dynamical degree $d_1(f)$ is larger than one or equivalently its topological and algebraic entropies are positive. Let $Q_n$ denote the set of isolated periodic points of period $n$ of $f$, counted with multiplicity. Then there is an invariant probability measure $\mu$ such that 
$$\lim_{n\to\infty} {1\over d_1(f)^n} \sum_{a\in Q_n} \delta_a =\mu.$$
\end{theorem}

The upper bound for the cardinality of $Q_n$ is given in Theorem \ref{th_d_k_per}. For the construction of a good family of periodic points, see Dujardin \cite{Dujardin}. The measure $\mu$ can be obtained as an intersection of invariant currents or as 
the intersection of the limit of $d_1(f)^{-n}[\Gamma_n]$ with the diagonal $\Delta$ of $X\times X$. We can also replace $Q_n$ in the theorem by the subset of saddle periodic points.

The reader will find more results for surfaces in \cite{DDG3, Dujardin, JonssonReschke} which are too technical to state in this paper. Notice that the upper bound for the number of isolated periodic points is crucial for the equidistribution property. It is sometimes overlooked in literature. We think that the strategy described in the introduction can be extended to holomorphic automorphisms or correspondences on projective manifolds of arbitrary dimension, whose actions on cohomology are simple. The needed techniques were recently developed in order to prove the following result.

Let $f$ be a polynomial automorphism of $\C^k$. We extend it to a birational map on the projective space $\P^k$ which is the natural compactification of $\C^k$. Denote by $I(f)$ and $I(f^{-1})$ the indeterminacy sets of $f$ and $f^{-1}$ respectively. They are analytic subsets of the hyperplane at infinity $\P^k\setminus \C^k$.
The following notion was introduced by the second author under the name of regular automorphisms, see \cite{Sibony}.

\begin{definition} \rm
We say that $f$ is a {\it H\'enon-type automorphism} if $f$ is not an automorphism of $\P^k$ and 
$I(f)\cap I(f^{-1})=\varnothing.$
\end{definition} 

This is a large family of of totally algebraically stable maps. The condition in the definition is easy to check. Moreover, in dimension 2, all polynomial automorphisms of $\C^2$ are conjugated to H\'enon-type maps or to elementary maps whose dynamics is simple to study, see Friedland-Milnor \cite{FM}. Consider a H\'enon-type map $f$ as above. It is not difficult to show that there is an integer $p$ such that $\dim I(f)=k-p-1$ and $\dim I(f^{-1})=p-1$. 
The action of  $f$  on cohomology is simple and $d_p(f)$ is the largest dynamical degree. Moreover, one can prove that $I(f^n)=I(f)$,  $I(f^{-n})=I(f^{-1})$ for $n\geq 1$, and then deduce that all periodic points of order $n$ of $f$ are isolated. 

The following result was obtained by Bedford-Lyubich-Smillie \cite{BLS2} for $k=2$ and by the authors in the general case \cite{DS6}.

\begin{theorem}[Bedford-Lyubich-Smillie, Dinh-Sibony] \label{th_per_Henon}
Let $f$ be  a H\'enon-type map on $\C^k$ as above and let $Q_n$ be the set of periodic points of period $n$ of $f$ counted with multiplicity. Then there is an invariant probability measure $\mu$ with compact support in $\C^k$ such that 
$$\lim_{n\to\infty} {1\over d_p(f)^n} \sum_{a\in Q_n} \delta_a =\mu.$$
\end{theorem}

This is the most difficult equidistribution property we obtained. 
In this theorem, we can also replace $Q_n$ by the subset of saddle periodic points. 
One may have periodic points at infinity but they are negligible for the above equidistribution property. Most of periodic points are in
a compact subset of $\C^k$.
The measure $\mu$ can be obtained as the intersection of invariant positive closed currents or as the intersection of the limit of $d_p(f)^{-n}[\Gamma_n]$ with the diagonal $\Delta$ of $\P^k\times \P^k$. It is known that the measure $\mu$ here is moderate. In particular, the support of $\mu$ is Zariski dense in $\C^k$.

In higher dimensions, the proof of the theorem requires a use of positive closed currents of arbitrary bidegree and is much more subtle than a simple technical issue. We follow
the strategy described in the introduction and the most difficult step is to show that 
$$\lim_{n\to\infty} \big(d_p(f)^{-n} (F^{n/2})^*[\Delta]\wedge[\Delta]\big) = \big(\lim_{n\to\infty} d_p(f)^{-n} (F^{n/2})^*[\Delta]\big)\wedge[\Delta].$$

In general, the two operations "taking a limit of currents" and "taking the wedge-product with a current"  do not commute. The main difficulty is the critical values of the projection from $\Gamma_n$ (to an imaginary space) following the direction of $\Delta$. In order to deal with this difficulty, we lift the dynamical system $F: \P^k\times\P^k\to \P^k\times\P^k$ to a suitable jet bundle associated to $\P^k\times \P^k$. This allows us to include the above "critical values" into the system. 
The behaviour of the lifts of $\Gamma_n$, or equivalently, of the orbit of $\Delta$ by $F$, to the jet bundle plays a crucial role in the proof of the above identity. A key point is that the volume of the lift of $\Gamma_n$ is of the same order of magnitude than the volume of $\Gamma_n$. 
One of the difficulties is that
the ergodicity properties we can get for the new system on the jet bundle are weaker than what we can obtain for $F$. To overcome these difficulties we use in particular the theory of densities of currents and some results by de Th\'elin \cite{deThelin2} on estimates of Lyapounov exponents in order to prove the above identity.

For the rest of this section, we will discuss the case of modular correspondences, where Problem 2 has a satisfactory solution. 
Let $G$ be a connected Lie group and let $\Lambda$ be a lattice in $G$. Define $\widehat X:=\Lambda\backslash G$. Let $g\in G$ be an element such that $g^{-1}\Lambda g$ is commensurable with $\Lambda$, that is, $\Lambda_g:=g^{-1}\Lambda g\cap \Lambda$ has finite index in $\Lambda$. Denote by $d_g$ this index.

The map $x\mapsto (x,gx)$ induces a map from $\Lambda_g\backslash G$ to $\widehat X\times\widehat X$. Let $\widehat \Gamma_g$ be its image. The natural projections from $\widehat \Gamma_g$ onto the factors of $\widehat X\times\widehat X$ define two unramified coverings of degree $d_g$. Therefore, it is the graph of a correspondence $\hat f$ on $\widehat X$ which is called an
 {\it (irreducible) modular} correspondence. When the group generated by $g$ and $\Gamma$ is dense in $G$, we say that $\hat f$ is an exterior correspondence. Let $K$ be a maximal compact subgroup of $G$. Since the left-multiplication and right-multiplication on $G$ commute, the correspondence $\hat f$ can descend to a correspondence $f$ on the locally symmetric space $X:=\widehat X/K$. We also say that $f$ is a modular correspondence on $X$. Several characterizations of modular correspondences on Hermitian locally symmetric spaces and related problems were considered in Clozel-Ullmo \cite{ClozelUllmo}, Huang-Yuan \cite{HY}, Mok \cite{Mok1, Mok2} and Mok-Ng \cite{MokNg}.  The following result was obtained in \cite{Dinh2}.

\begin{theorem}[Dinh] \label{th_per_modular}
Let $f$ be an exterior modular correspondence on a locally symmetric space $X$ as above. Let $Q_n$ denote the set of isolated periodic points of period $n$ of $f$. Assume also that $f$ has at least one isolated periodic point. Then 
$$\lim_{n\to\infty} {1\over \#Q_n} \sum_{a\in Q_n} \delta_a =\mu,$$
where $\mu$ is the Haar (probability) measure on $X$. 
\end{theorem}

Note that the theorem holds also when $X$ is not Hermitian or not compact. The proof of this theorem uses an analysis of the inverse branches of $f^n$ on balls of $X$, similar to the one in Theorem \ref{th_eq_dominant}. However, since $f$ is locally isometric with respect to an invariant Riemannian metric on $X$, the problem is technically simpler. A key point in the proof is the equidistribution of the orbits of points with respect to the measure $\mu$ which was obtained by Clozel-Otal in \cite{ClozelOtal}.

\section{Equidistribution of orbits of points and varieties} \label{section_varieties}

As mentioned in the introduction, the equidistribution of periodic points is closely related to the equidistribution of the orbits of points and varieties. In this section, we will discuss some situations where we have a satisfactory answer to the following problem.
Let $f$ be a dominant meromorphic correspondence on a compact K\"ahler manifold $X$ of dimension $k$ as above. 

\begin{problem}
Let $V$ be an analytic subset of $X$ of pure codimension $p$. Study the convergence and rate of convergence of the currents of integration on $f^{-n}(V)$, properly normalized, when $n$ goes to infinity.
\end{problem}

Note that for convenience, we state the problem for the negative orbit of $V$ by $f$. The study of the positive orbit of $V$ by $f$  is equivalent to the study of the negative orbit of $V$ by $f^{-1}$. 

Recall that the mass of a positive closed current depends only on its cohomology class. So the mass of the current $[f^{-n}(V)]$ grows at most as the norm of the operator $(f^n)^*$ on $H^{p,p}(X,\C)$. Hence, concretely, we want to study the convergence and the rate of convergence of $d_n^{-1}[f^{-n}(V)]$ for a suitable sequence $(d_n)$ of positive numbers which is equivalent to the sequence  of norms $\|(f^n)^*:H^{p,p}(X,\C)\to H^{p,p}(X,\C)\|$. 

In the situations we will consider, the correspondence is algebraically $p$-stable. We also have  $d_p>d_l$ for $l<p$ except for the case of modular correspondences. Unless we use strong geometric hypotheses or other conditions, the case without algebraic stability seems not to be accessible for the moment. Without the condition $d_p>d_l$, the nature of the problem changes. The reader can, for example, compare the case of the orbits of points by a H\'enon-type automorphism with the statements below. We first discuss the cases of orbits of points and hypersurfaces which only require classical tools from complex analysis and geometry for the proofs.

\begin{theorem}
Let $f$ be a non-invertible holomorphic map on $\P^k$. Let $d_k(f)$ denote its topological degree. Then
\begin{enumerate}
\item There is a maximal proper analytic subset $\Ec$ of $\P^k$ which is totally invariant, i.e. we have  $f^{-1}(\Ec)=f(\Ec)=\Ec$.
\item There is an invariant moderate probability measure $\mu$ of $f$ (called equilibrium measure) such that 
$$\lim_{n\to\infty} {1\over d_k(f)^n} \sum_{a\in f^{-n}(z)} \delta_a =\mu \quad \text{if and only if} \quad z\not\in \Ec,$$
 where the points in $f^{-n}(z)$ are counted with multiplicity.
 \item The convergence in the previous statement for $z\not\in \Ec$ is exponentially fast. More precisely, if $\varphi$ is a H\"older continuous function on $\P^k$ then the integral
 $$\Big \langle {1\over d_k(f)^n} \sum_{a\in f^{-n}(z)} \delta_a - \mu, \varphi\Big\rangle$$ 
 tends to $0$ exponentially fast as $n$ goes to infinity.
\end{enumerate}
\end{theorem}

This result was obtained by the authors using strong compactness properties of quasi-p.s.h. functions (see Theorem \ref{th_Skoda}) which are used as test functions for measures and also some analysis of the critical set in order to localize the exceptional set $\Ec$. We refer the reader to \cite{DS7} for more precise statements and history of the problem. The similar statement with the same proof holds  for algebraically expanding holomorphic maps on compact K\"ahler manifolds.
The convergence in the second assertion has been obtained by Forn\ae ss and the second author for $z$ outside a pluripolar set and by Briend-Duval for $z$ outside a countable union of analytic subsets, see \cite{BD2, FS2}. For the case of dimension 1, see  
Brolin, Freire-Lopes-Ma{\~n}{\'e} and Lyubich \cite{Brolin, FLM, Lyubich}.

Consider now a more general setting. Let $f$ be an algebraically expanding map on a compact K\"ahler manifold $X$ of dimension $k$. Let $d_k(f)$ denote its topological degree. Let $I(f)$ and $I(f^{-1})$ denote the indeterminacy sets of $f$ and $f^{-1}$. 
Denote by $I_\infty(f)$ the positive orbit of $I(f)$ and $I_\infty(f^{-1})$ the one of $f^{-1}$. We will consider points outside those sets. The following result was obtained by Nguyen, Truong and the first author using a geometrical approach, originally used in dimension 1, see \cite{DNT2}. The convergence has been obtained by the authors for $z$ outside a pluripolar set and by Guedj for $z$ outside a countable union of analytic sets, see \cite{DS8, Guedj1}.

\begin{theorem}
Let $f$ be an algebraically expanding map as above. There is an invariant probability measure $\mu$ of $f$ (equilibrium measure) and a proper analytic subset $\Ec$ of $X$ such that if $z$ is outside $I_\infty(f)\cup I_\infty(f^{-1})$, then
$$\lim_{n\to\infty} {1\over d_k(f)^n} \sum_{a\in f^{-n}(z)} \delta_a =\mu \quad \text{if and only if} \quad z\not\in \Ec,$$
 where the points in $f^{-n}(z)$ are counted with multiplicity.
\end{theorem}

Note that when $z$ is outside $I_\infty(f)\cup I_\infty(f^{-1})$, the condition  $z\not\in \Ec$ is equivalent to the condition that $z$ doesn't belong to the orbit of $\Ec$ because $\Ec$ satisfies a strong invariant property : $f^{-1}(\Ec\setminus I(f^{-1}))\subset \Ec$. The result can be extended to algebraically expanding correspondences but in this case $\Ec$ doesn't satisfy such an invariant property and therefore we need to replace the condition $z\not\in \Ec$ with the condition that $z$ is not in the positive orbit of $\Ec$ by $f$.

Consider now the case of a modular correspondence which plays the central role in the proof of Theorem \ref{th_per_modular}. The following result was obtained by Clozel-Otal in \cite{ClozelOtal}, see also Clozel-Ullmo \cite{ClozelUllmo}.

\begin{theorem}[Clozel-Otal]
Let $f, X$ and $\mu$ be  as in Theorem \ref{th_per_modular}. Let $d$ be the topological degree of $f$. Then for every $z\in X$ we have
$$\lim_{n\to\infty} {1\over d^n} \sum_{a\in f^{-n}(z)} \delta_a =\mu.$$
\end{theorem}

We discuss now the case of orbits of hypersurfaces. The following result was obtained by Favre-Jonsson in the case of dimension 2 and by the authors in the general case \cite{DS9, FavreJonsson}. 

\begin{theorem}[Favre-Jonsson, Dinh-Sibony]
Let $f$ be a non-invertible holomorphic map on $\P^k$. Let $d_1(f)$ denote its first dynamical degree. 
Then there is a totally invariant proper analytic set $\Ec_0$ such that if $V$ is a hypersurface of degree $\deg(V)$ of $\P^k$ which doesn't contain any component of $\Ec_0$ then 
$$ \lim_{n\to\infty} {1\over \deg(V)d_1(f)^n} (f^n)^*[V]=T,$$
where $T$ is the dynamical Green $(1,1)$-current of $f$. 
\end{theorem}

So the negative orbit of $V$ is equidistributed with respect to the current $T$. There are different ways to construct this invariant positive closed $(1,1)$-current. The above theorem shows that it is the limit of $\deg(V)^{-1}d_1(f)^{-n} (f^n)^*[V]$ for generic hypersurfaces $V$. If $\alpha$ is a closed $(1,1)$-form with bounded coefficients in $\P^k$ which is cohomologous to a hyperplane, then $d_1(f)^{-n} (f^n)^*(\alpha)$ also converges to $T$ in the sense of currents, see e.g. \cite{DS10} for more details.

Note that for generic hypersurfaces $V$, Taflin proved that the convergence of the currents $\deg(V)^{-1}d_1(f)^{-n} (f^n)^*[V]$ towards $T$ is exponentially fast, see \cite{Taflin} for a precise statement. Note also that all the above results can be extended to H\'enon-type automorphisms of $\C^k$, see \cite{DS3, Taflin}. 

The key point in the proof is to use normalized quasi-potentials of positive closed $(1,1)$-currents. The convergence of these currents is equivalent to the convergence of their  normalized quasi-potentials. Compactness of quasi-p.s.h. functions and other classical tools from complex geometry allow to obtain the above results.

Consider now the case of analytic sets $V$ of arbitrary dimension. The theory of super-potentials described in Section \ref{section_current} substitutes the use of quasi-p.s.h. functions. In particular, we use Theorem \ref{th_super_pot} instead of Theorem \ref{th_Skoda}.  We however need to overcome more technical difficulties than in the  hypersurfaces case.

Let $\Hc_d$ denote the family of all holomorphic self-maps of $\P^k$ such that the first dynamical degree is $d>1$. This can be identified to a Zariski open subset of a projective space. The following result was obtained by the authors in \cite{DS3}, see also Ahn \cite{Ahn} for some extension.

\begin{theorem}[Dinh-Sibony] \label{th_equi_variety}
There is an explicit dense Zariski open subset $\Hc_d'$ of $\Hc_d$ such that for every $f$ in $\Hc_d'$ and every analytic subset $V$ of pure codimension $p$ of $\P^k$ we have
$$\lim_{n\to\infty} {1\over d^p\deg(V)} (f^n)^*[V]=T^p,$$
where $T^p$ is the $p$-th power of the  dynamical Green $(1,1)$-current of $f$. Moreover, the convergence is exponentially fast and uniform on $V$.  
\end{theorem}

In the case of H\'enon-type maps, the critical values of $f^n$ are easier to understand and we obtain in the same way a stronger property. Let $f$ be a H\'enon-type automorphism of $\C^k$ as in Theorem \ref{th_per_Henon}. Recall that $I(f^{-1})$ is attractive for $f$. Let $\Uc(f)$ denote the basin of $I(f^{-1})$ which is an open neighbourhood of $I(f^{-1})$ in $\P^k$. The set $\Kc(f):=\C^k\setminus \Uc(f)$ is the set of all points $z\in\C^k$ whose positive orbits by $f$ are bounded in $\C^k$. It is also known that the closure $\overline {\Kc(f)}$ of $\Kc(f)$ in $\P^k$ is the union of $\Kc(f)$ with $I(f)$. Recall that $\dim I(f)=k-p-1$ and $\dim I(f^{-1})=p-1$.
The following result was obtained by Forn\ae ss and the second author for the case of dimension 2 in \cite{FS2} and by the authors for the general case in \cite{DS3}, see also \cite{BLS1, FS3}.

\begin{theorem}[Forn\ae ss-Sibony, Dinh-Sibony] \label{th_Henon_rigid}
Let $V$ be an analytic subset of pure dimension $k-p$ and degree $\deg(V)$ in $\P^k$ such that $V\cap I(f^{-1})=\varnothing$. 
Then $\deg(f)^{-1}d_p(f)^{-n} (f^n)^*[V]$ converges exponentially fast to a positive closed $(p,p)$-current $T(f)$. Moreover, $T(f)$ is the unique positive closed $(p,p)$-current of mass $1$ with support in $\overline{\Kc(f)}$. 
\end{theorem}

In the case of dimension 2, the last statement also holds if we replace positive closed currents by positive $\ddc$-closed currents, see \cite{DS11}.
Note that we can also apply the theorem for $f^{-1}$. 
The equilibrium measure $\mu$ considered in Theorem \ref{th_per_Henon} can be obtained as the intersection of $T(f)$ and $T(f^{-1})$. 
In general, the map $F=(f,f^{-1})$ on $\C^k\times \C^k$ is not a H\'enon-type map but its dynamics can be studied in a similar way. A version of the last theorem for $F$ is a crucial point in the proof of Theorem \ref{th_per_Henon} following the general strategy described in the introduction. In particular, the sequence of currents $d_p(f)^{-n}[\Gamma_n]$ converges to the current $T(f)\otimes T(f^{-1})$. The intersection of the last current with the diagonal of $\P^k\times\P^k$ can be identified with the measure $\mu$.

Note that one can obtain similar properties for holomorphic automorphisms of a compact K\"ahler manifold, see \cite{Cantat, DS4, DS11} for details.

\section{Open problems related to equidistribution} \label{section_problem}

In previous sections, we already mentioned several general open problems. We will give in this section some other concrete questions which require new ideas or important technical tools. In order to keep the presentation simple, some problems are stated for particular families of maps but the reader can easily generalize them to other maps and correspondences or their restrictions to certain open subsets, e.g. to the basin of an attracting set, see \cite{Dinh3, Taflin2}.

\begin{problem} \label{prob_per_aut}
Let $f$ be a holomorphic automorphism of positive entropy on a compact K\"ahler surface or more generally a holomorphic automorphism of positive entropy on a compact K\"ahler manifold with a simple action on cohomology. Show that its periodic points are equidistributed with respect to an invariant probability measure.
\end{problem} 

We have seen in Theorem \ref{th_per_surface} that this property holds when $X$ is a projective surface. Likely, the techniques in the proof of Theorem \ref{th_per_Henon} allow to solve the problem  when $X$ is a projective manifold. However, the general case of K\"ahler manifolds requires a new technical idea whose interest may be greater than the solution of the last question, see e.g. the next problem.  
 
Consider a holomorphic endomorphism $f$ of $\P^k$ with the first dynamical degree $d>1$. Let $Q_n$ denote the set of periodic points of period $n$ of $f$ counted with multiplicity. Note that in this case, all the periodic points of period $n$ are isolated and we have $\#Q_n=d^{kn}+O(d^{(k-1)n})$. 

\begin{problem}
Study the rate of the following convergence
$$\lim_{n\to\infty} {1\over d^n}\sum_{a\in Q_n}\delta_a=\mu$$
given in Theorem \ref{th_eq_dominant}.
\end{problem}

The same ideas used in the proof of Theorem  \ref{th_per_Henon} can be applied to get the last convergence property. The situation is even much simpler in this case. Several steps of the proof are already quantitative but at the end, we use some 
non-quantitative arguments from complex geometry. One should improve the techniques in order to get a quantitative result. This part may be related to Problem \ref{prob_per_aut}. 

Recall that $T^p$ is the  dynamical Green $(p,p)$-current of $f$, see Theorem \ref{th_equi_variety}. Its support is called {\it the Julia set of order $p$} of $f$. De Th\'elin and the first author proved that the topological entropy of $f$ on any compact subset disjoint from $\supp(T^p)$ is at most equal to $(p-1)\log d$, see \cite{deThelin3, Dinh3}. The following question is directly related to the last one.
 
\begin{problem}
Let $K$ be a compact set which is disjoint from $\supp(T^p)$. Do we always have 
$$\limsup_{n\to\infty} {1\over n}\log \# (Q_n\cap K)\leq (p-1)\log d.$$
\end{problem} 

Note that Forn\ae ss and the second author showed that one may have an infinite number of periodic points outside the support of the equilibrium measure $\mu=T^k$ when $k\geq 2$, see \cite{FS4}. It is well-known that in dimension $k=1$ there are only finitely many of such points.

The following folklore conjecture has been solved in some particular  cases, see Hwang-Nakayama \cite{HN} and the references therein. We believe that the existence of repelling periodic points can be used to study the problem.
Indeed, the existence of repelling periodic points obtained in Theorem \ref{th_eq_dominant} implies that $X$ contains infinitely many non-degenerate holomorphic images of $\C^n$, where $\dim X=n$.

\begin{conjecture}
Let $X$ be a Fano manifold with Picard number $1$. Assume that $X$ admits a non-invertible holomorphic endomorphism. Then $X$ should be a projective space.
\end{conjecture}

Concerning the equidistribution property in Theorem \ref{th_equi_variety}, the following conjecture was stated in \cite{DS9}. It is open for $2\leq p\leq k-1$.

\begin{conjecture}
Let $V$ be an irreducible analytic subset of dimension $k-p$ of $\P^k$. Assume that $V$ is generic in the sense that 
$V\cap E=\varnothing$ or $\codim V\cap E=p+\codim E$ 
for any irreducible component $E$ of every totally invariant analytic
subset of $\P^k$.
Then 
$$\lim_{n\to\infty} d^{-pn}(f^n)^*[V]=T^p.$$
\end{conjecture}

We can also investigate the rate of the last convergence. Note that the authors proved that there are only finitely many of analytic subsets $E$ of $\P^k$ which are totally invariant, i.e. $f^{-1}(E)=f(E)=E$, see e.g. \cite{DS10}. In comparison with the known case
of hypersurface with $p=1$, super-potentials can replace quasi-p.s.h. functions. However, we still need to extend some arguments from classical complex geometry (e.g. Lojasiewicz's inequality) to the space of positive closed currents which is of infinite dimension.

Note that we can extend the conjecture to the case of H\'enon-type maps by considering totally invariant sets for the restriction of the map to the indeterminacy set of its inverse.

Properties of totally invariant analytic subsets of $\P^k$ may be useful for the last question. The following problem is still open in dimension $k\geq 3$. For the case of dimension 2, see \cite{CLN, FS0} and also \cite{AC, Zhang}.

\begin{problem}
Let $E$ be an analytic subset of $\P^k$ which is totally invariant by $f$. Is $E$ always a union of linear analytic subspaces of $\P^k$ ?
Find an upper bound for the degree of $E$ (or for $\#E$ when $E$ is finite) which depends only on $k$. Note that we don't assume that $E$ is irreducible.
\end{problem}

We will end this section with the particular case of H\'enon-type maps. Let $f$ be a H\'enon-type map as in Theorem \ref{th_per_Henon}. Theorem \ref{th_Henon_rigid} shows that $\overline{\Kc(f)}$ satisfies a very strong rigidity property.
For $p=k-1$, i.e. $\dim I(f)=0$, we know that there are holomorphic entire maps with values in $\overline{\Kc(f)}$. Using Nevanlinna's theory, one can produce positive closed current of bidimension $(1,1)$ with support in $\overline{\Kc(f)}$. By Theorem \ref{th_Henon_rigid}, this current is proportional to the dynamical Green $(k-1,k-1)$-current of $f$ whose support is contained the boundary of $\Kc(f)$. This support is exactly   the boundary of $\Kc(f)$ when $k=2$ and $p=1$.
In the following problem, the first question is open for $p<k-1$, i.e. $\dim I(f)>0$,  and the second one is open for $p>1$.

\begin{problem}
Let $\tau:\C^{k-p}\to \Kc(f)$ be a non-degenerate holomorphic map. Does the closure of $\tau(\C^{k-p})$ contain the support of the dynamical Green $(k-p,k-p)$-current of $f$ ? Is the support of this current equal to the boundary of $\Kc(f)$ ?
\end{problem}

Nevanlinna's theory is well developed for holomorphic maps on $\C$ but the theory doesn't seem to work for general holomorphic maps on $\C^{k-p}$. Some aspects can be extended to higher dimension under conditions which are difficult to check. The last problem may lead to a setting where one can develop the theory with applications.


\small

\begin{thebibliography}{99}

\bibitem{Ahn}
Ahn T., 
Equidistribution in higher codimension for holomorphic endomorphisms of $\P^k$. 
{\it Trans. Amer. Math. Soc.} {\bf 368} (2016), no. 5, 3359-3388. 

\bibitem{AC}
Amerik E., Campana F.,
Exceptional points of an endomorphism of the projective plane. 
{\it Math. Z.} {\bf 249} (2005), no. 4, 741-754. 



\bibitem{BLS2}
Bedford E., Lyubich M., Smillie J., 
Distribution of periodic points of polynomial diffeomorphisms of $\bold C^2$. {\it Invent. Math.} {\bf 114} (1993), no. 2, 277-288.


\bibitem{BLS1}
Bedford E., Lyubich M., Smillie J., 
Polynomial diffeomorphisms of $\C^2$. IV. The measure of maximal entropy and laminar currents. {\it Invent. Math.} {\bf 112} (1993), no. 1, 77-125.

\bibitem{BBD}
Berteloot F., Bianchi F., Dupont C.,
Dynamical stability and Lyapunov exponents for holomorphic endomorphisms of $\P^k$. {\it Preprint} (2014). {\tt arXiv:1403.7603}

\bibitem{BD1} 
 Briend J.-Y., Duval J., Exposants de Liapounoff et distribution des points p\'eriodiques d'un endomorphisme de $\mathbb{C}\mathbb{P}^k$.  {\it  Acta Math.}  {\bf 182} (1999), no. 2, 143-157.
 
 \bibitem{BD2}
Briend J.-Y., Duval J., Deux caract{\'e}risations de la mesure
d'{\'e}quilibre d'un endomorphisme de ${\rm P}\sp k(\bold C)$. {\it
Publ. Math. Inst. Hautes {\'E}tudes Sci.} {\bf 93} (2001), 145-159. Erratum : {\it
Publ. Math. Inst. Hautes {\'E}tudes Sci.} {\bf 109} (2009), 295-296.


\bibitem{Brolin}
Brolin H., Invariant sets under iteration of rational functions.
{\it Ark. Mat.} {\bf 6} (1965), 103-144.

\bibitem{Cantat}
Cantat S., Dynamique des automorphismes des surfaces K3. {\it Acta Math.} {\bf 187} (2001), 1-57.

\bibitem{CLN}
Cerveau D., Lins Neto A., Hypersurfaces exceptionnelles des endomorphismes de $\C\P(n)$. {\it Bol. Soc. Brasil. Mat. (N.S.)} {\bf 31} (2000), no. 2, 155-161.


\bibitem{ClozelOtal}
Clozel L., Otal J.-P.,
Unique ergodicit\'e des correspondances modulaires. {\it Essays on geometry and related topics}, Vol. 1, 2, 205-216, 
Monogr. Enseign. Math. {\bf 38}, {\it Enseignement Math.}, Geneva, 2001.

\bibitem{ClozelUllmo}
Clozel L., Ullmo E., Correspondances modulaires et mesures invariantes. {\it J. Reine Angew. Math.} {\bf 558} (2003), 47-83.


\bibitem{Demailly}
Demailly J.-P., {\it Complex analytic and differential geometry}. Available at \\
{\tt www.fourier.ujf-grenoble.fr/$\sim$demailly}. 


\bibitem{deThelin3}
De Th{\'e}lin H.,  Sur la construction de mesures selles. {\it
  Ann. Inst. Fourier}  {\bf 56}  (2006),  no. 2, 337-372.
  
\bibitem{deThelin2} 
De Th{\'e}lin H.,  
Sur les exposants de Lyapounov des applications m\'eromorphes. 
{\it Invent. Math.} {\bf 172}  (2008), no. 1, 89-116.

  
\bibitem{DDG3}
 Diller J., Dujardin R., Guedj V., Dynamics of meromorphic maps with small topological degree III: geometric currents and ergodic theory. {\it Ann. Sci. \'Ec. Norm. Sup\'er. (4)} {\bf 43} (2010), no. 2, 235-278.
 
 \bibitem{Dinh1}
Dinh T.-C., Suites d'applications m{\'e}romorphes multivalu{\'e}es et courants
laminaires.  {\it J. Geom. Anal.} {\bf  15}  (2005),  207-227.

\bibitem{Dinh3}
Dinh T.-C., Attracting current and equilibrium measure for attractors
on $\Bbb P\sp k$. {\it  J. Geom. Anal.} {\bf  17}  (2007),  no. 2, 227-244.

 \bibitem{Dinh2}
Dinh T.-C., 
Equidistribution of periodic points for modular correspondences. 
{\it J. Geom. Anal.} {\bf 23} (2013), no. 3, 1189-1195. 

\bibitem{DN1}
Dinh T.-C., Nguyen V.-A., The mixed Hodge-Riemann bilinear relations for compact K\"ahler manifolds. {\it Geom. Funct. Anal.} {\bf 16} (2006), no. 4, 838-849.

\bibitem{DNS1}
Dinh T.-C., Nguyen V.-A., Sibony N., Exponential estimates for
plurisubharmonic functions and stochastic dynamics. {\it J. Diff. Geom.} 
{\bf 84} (2010), no. 3, 465-488.

\bibitem{DNT1}
Dinh T.-C., Nguyen V.-A., Truong T. T.,
On the dynamical degrees of meromorphic maps preserving a fibration. {\it Comm. Contemporary. Math.}
{\bf 14}, No. 6 (2012). DOI: 10.1142/S0219199712500423

\bibitem{DNT2}
Dinh T.-C., Nguyen V.-A., Truong T. T., Equidistribution for meromorphic maps with dominant topological degree. {\it Indiana J. Math.} {\bf 64} (2015), no. 6, 1805-1828. 

\bibitem{DNT3}
Dinh T.-C., Nguyen V.-A., Truong T. T., 
Growth of the number of periodic points for meromorphic maps. {\it Preprint}, 2016. {\tt arXiv:1601.03910}


  \bibitem{DS0}   Dinh T.-C., Sibony N.,
  Dynamique des applications d'allure polynomiale. {\it J. Math. Pures Appl. (9)}  {\bf  82} (2003), no. 4, 367-423. 

\bibitem{DS1}
Dinh T.-C., Sibony N., Regularization of currents and entropy. {\it Ann. Sci. Ecole
Norm. Sup.}  {\bf 37} (2004), 959-971.

\bibitem{DS13}
Dinh T.-C., Sibony N., Groupes commutatifs d'automorphismes
holomorphes d'une vari{\'e}t{\'e} k{\"a}hl{\'e}rienne compacte.
{\it Duke Math. J.} {\bf 123} (2004),  311-328.

\bibitem{DS8}
Dinh T.C., Sibony N., Distribution des valeurs d'une suite de
transformations m{\'e}romorphes et applications. 
{\it  Comment. Math. Helv.} {\bf 81} (2006), 221-258.

\bibitem{DS2}
Dinh T.-C., Sibony N.,  Pull-back currents by holomorphic maps. {\it Manuscripta Math.} {\bf 123} (2007), no. 3, 357-371.

\bibitem{DS12}
Dinh T.-C., Sibony N., Upper bound for the topological entropy of a
meromorphic correspondence.  {\it  Israel J. Math.} {\bf  163}  (2008), 29-44.

\bibitem{DS9}
Dinh T.-C., Sibony N., Equidistribution towards the Green current for
holomorphic maps. {\it Annales Sci. ENS} {\bf 41} (2008), 307-336.

\bibitem{DS3}
Dinh T.-C., Sibony N., Super-potentials of positive closed currents, intersection theory
and dynamics. {\it Acta Math.} {\bf 203} (2009), no. 1, 1-82.

\bibitem{DS7}
Dinh T.-C., Sibony N., Equidistribution speed for endomorphisms of
projective spaces.  {\it Math. Ann.} {\bf 347} (2010), no. 3, 613-626.

\bibitem{DS10}
Dinh T.-C., Sibony N., Dynamics in several complex variables:
endomorphisms of projective spaces and polynomial-like mappings,
{\it Holomorphic dynamical systems}, 165-294, Lecture Notes in Math., {\bf 1998}, Springer, Berlin, 2010.

\bibitem{DS4}
Dinh T.-C., Sibony N., 
Super-potentials for currents on compact K\"ahler manifolds and dynamics of automorphisms. {\it J. Algebraic Geom.} {\bf 19} (2010), no. 3, 473-529.

\bibitem{DS5}
Dinh T.-C., Sibony N.,  Density of positive closed currents, a theory of non-generic intersections. {\it Preprint}, new version 2015. {\tt arXiv:1203.5810}

\bibitem{DS6}
Dinh T.-C., Sibony N., 
Equidistribution of saddle periodic points for H\'enon-type automorphisms of $\C^k$. {\it Math. Ann.} {\bf 366} (2016), no. 3-4, 1207-1251.


\bibitem{DS11}
Dinh T.-C., Sibony N., 
Rigidity of Julia sets for H\'enon type maps. {\it J. Mod. Dyn.} {\bf 8} (2014), no. 3-4, 499-548. 


\bibitem{Dujardin}
Dujardin R., Laminar currents and birational dynamics. {\it Duke Math. J.} {\bf 131}, no. 2 (2006), 219-247. 


\bibitem{EOY}
Esnault H., Oguiso K.,  Yu X., Automorphisms of elliptic K3 surfaces and Salem numbers of maximal degree. {\it Preprint} (2014). {\tt arXiv:1411.0769}

\bibitem{ES}
Esnault H.,  Srinivas V., Algebraic versus topological entropy for surfaces over finite fields. {\it Osaka J. Math.} {\bf 50} (2013), no 3, 827-846.


 \bibitem{Favre} 
 Favre C.,  Points p\'eriodiques d'applications birationnelles de $\P^2.$  {\it Ann. Inst. Fourier (Grenoble)} {\bf  48} (1998), no. 4, 999-1023.
 
 \bibitem{FavreJonsson}
Favre C., Jonsson M., Brolin's theorem for curves in two complex
dimensions. {\it  Ann. Inst. Fourier} {\bf 53}  (2003), no. 5,
1461-1501.

\bibitem{FS3}
Forn\ae ss J.-E., Sibony N.,  Complex H\'enon mappings in $\C^2$ and Fatou-Bieberbach domains. {\it Duke Math. J.} {\bf 65} (1992), no. 2, 345-380

\bibitem{FS0}
Forn\ae ss J.-E., Sibony N., 
 Complex dynamics in higher dimension. I. In: Complex Analytic Methods in Dynamical Systems, Rio de Janeiro, 1992, {\it Ast\'erisque}, vol. {\bf 222}, Soc. Math. France, 1994, pp. 201-231.
 
\bibitem{FS2}
Forn\ae ss J.-E., Sibony N., 
Complex dynamics in higher dimensions. Notes partially written by Estela A. Gavosto. {\it NATO Adv. Sci. Inst. Ser. C Math. Phys. Sci.}, {\bf 439}, Complex potential theory (Montreal, PQ, 1993), 131-186, Kluwer Acad. Publ., Dordrecht, 1994.
 
\bibitem{FS1}
Forn\ae ss J.-E., Sibony N., Complex dynamics in higher dimension. II. {\it Modern methods in complex analysis (Princeton, NJ, 1992)}, 135-182, {\it Ann. of Math. Stud.} {\bf 137}, Princeton Univ. Press, Princeton, NJ, 1995.

\bibitem{FS4}
Forn\ae ss J.-E., Sibony N.,
Dynamics of $\P^2$
(examples). Laminations and foliations in dynamics, geometry and topology (Stony Brook, NY, 1998), 47-85, {\it Contemp. Math.}, {\bf 269}, Amer. Math. Soc., Providence, RI, 2001.

\bibitem{FLM}
Freire A., Lopes A., Ma{\~n}{\'e} R., An invariant measure for
rational maps. {\it Bol. Soc. Brasil. Mat.} {\bf 14} (1983), no.
1, 45-62.

\bibitem{FM}
Friedland S., Milnor J.,
Dynamical properties of plane polynomial automorphisms.
{\it Ergodic Theory Dynam. Systems} {\bf 9} (1989), no. 1, 67-99. 

\bibitem{Gromov1}
Gromov M., Convex sets and K{\"a}hler manifolds.
\textit{Advances in differential geometry and topology}.
Word Sci. Publishing, Teaneck, NJ, 1998, 1-38.

\bibitem{Gromov2}
Gromov M.,  On the entropy of holomorphic maps.
\textit{Enseignement Math.} \textbf{49} (2003), 217-235. {\it Manuscript} (1977).

\bibitem{Guedj1}
Guedj V.,
Ergodic properties of rational mappings with large topological degree. {\it Ann. of Math.} (2) {\bf 161} (2005), no. 3, 1589-1607.

\bibitem{Hormander}
H{\"o}rmander L., {\it An introduction to complex analysis in several
  variables}. Third edition, 
North-Holland Mathematical Library, {\bf 7}, North-Holland Publishing Co., Amsterdam, 1990. 

\bibitem{HY}
Huang X., Yuan Y.,
Holomorphic isometry from a K\"ahler manifold into a product of complex projective manifolds. 
{\it Geom. Funct. Anal.} {\bf 24} (2014), no. 3, 854-886. 


\bibitem{HN}
Hwang J.-M., Nakayama N.,
On endomorphisms of Fano manifolds of Picard number one. 
{\it Pure Appl. Math. Q.} {\bf 7} (2011), no. 4, Special Issue: In memory of Eckart Viehweg, 1407-1426. 

\bibitem{IwasakiUehara}  
 Iwasaki K., Uehara T., Periodic points for area-preserving birational maps of surfaces. 
{\it Math. Z.}  {\bf 266} (2010), no. 2, 289-318.  

\bibitem{JonssonReschke}
Jonsson M., Reschke J.,
On the complex dynamics of birational surface maps defined over number fields. {\it  J. Reine Angew. Math.}
(to  appear), {\tt  arXiv:1505.03559}.

\bibitem{Kaloshin}
 Kaloshin V. Yu.
Generic diffeomorphisms with superexponential growth of number of periodic orbits. 
{\it Comm. Math. Phys.} {\bf 211} (2000), no. 1, 253-271.

\bibitem{Kaufmann}
Kaufmann L., A Skoda-type integrability theorem for singular Monge-Amp\`ere measures. {\it Michigan Mathematical Journal}, to appear.


\bibitem{Lyubich}
Lyubich M. Ju., Entropy properties of rational endomorphisms of
the Riemann sphere. {\it Ergodic Theory Dynam. Systems} {\bf  3}
(1983), no. 3, 351-385.

\bibitem{Mok1}
Mok N., 
Local holomorphic isometric embeddings arising from correspondences in the rank-1 case. {\it Contemporary trends in algebraic geometry and algebraic topology (Tianjin, 2000)}, 155-165, 
Nankai Tracts Math., 5, {\it World Sci. Publ., River Edge, NJ,} 2002.
 
\bibitem{Mok2}
Mok N.,  
 Extension of germs of holomorphic isometries up to normalizing constants with respect to the Bergman metric. 
{\it J. Eur. Math. Soc.} {\bf 14} (2012), no. 5, 1617-1656. 

 
\bibitem{MokNg}
Mok N.,  Ng S.-C.,  Germs of measure-preserving holomorphic maps from bounded symmetric domains to their Cartesian products. {\it J. Reine Angew. Math.} {\bf 669} (2012), 47-73.

\bibitem{Nguyen}
Nguyen V.-A., Green currents for quasi-algebraically stable meromorphic self-maps of $\P^k$. {\it Publ. Mat.} {\bf 56} (2012), no. 1, 127-146.
 
\bibitem{Oguiso1}
Oguiso K., 
Pisot units, Salem numbers and higher dimensional projective manifolds with primitive automorphisms of positive entropy. {\it Preprint} (2016). {\tt arXiv:1608.03122}

\bibitem{OT}
Oguiso K., Truong T.T.,  Explicit examples of rational and Calabi-Yau threefolds with primitive automorphisms of positive entropy. {\it Kodaira Centennial issue of the Journal of Mathematical Sciences}, the University of Tokyo, {\bf 22} (2015) 361-385.

\bibitem{RS}
Russakovskii A., Shiffman B., Value distribution for sequences of rational mappings and complex dynamics. {\it Indiana Univ. Math. J.} {\bf 46} (1997), no. 3, 897-932. 

 \bibitem{Saito}   
 Saito S.,  General fixed point formula for an algebraic surface and the theory of Swan representations for two-dimensional local rings. 
{\it Amer. J. Math.} {\bf 109} (1987), no. 6, 1009-1042.   

\bibitem{Sibony} Sibony N., Dynamique des applications rationnelles de $\mathbb{P}^k$. \textit{Panoramas et Synth{\`e}ses} \textbf{8} (1999), 97-185.



\bibitem{Siu}
Siu Y.T., Analyticity of sets associated to Lelong numbers and the extension of closed positive currents. {\it Invent. Math.}, {\bf  27}, (1974), 53-156.


\bibitem{Taflin}
Taflin J., 
Equidistribution speed towards the Green current for endomorphisms of $\P^k$.
{\it Adv. Math.} {\bf 227} (2011), no. 5, 2059-2081. 


\bibitem{Taflin2}
Taflin J., 
Speed of convergence towards attracting sets for endomorphisms of $\P^k$. {\it Indiana Univ. Math. J.} {\bf 62} (2013), no. 1, 33-44.

\bibitem{Truong1}
Truong T.T.,  (Relative) dynamical degrees of rational maps over an algebraic closed field. {\it Preprint} (2015). {\tt  arXiv:1501.01523}

\bibitem{Truong2}
Truong T.T.,  Relative dynamical degrees of correspondences over a field of arbitrary characteristic. {\it Preprint} (2016). {\tt  arXiv:1605.05049}

\bibitem{Truong3}
Truong T.T.,  Relations between dynamical degrees, Weil's Riemann hypothesis and the standard conjectures. {\it Preprint} (2016).
{\tt arXiv:1611.01124}

\bibitem{Voisin} 
Voisin C., {\it Th{\'e}orie de Hodge et g{\'e}om{\'e}trie alg{\'e}brique complexe}. Cours Sp{\'e}cialis{\'e}s {\bf 10}, Soci{\'e}t{\'e} Math{\'e}matique de France, Paris, 2002. 

\bibitem{Vu2}
Vu D.-V.,
Intersection of positive closed currents of higher bidegree. {\it Preprint} (2015). {\tt arXiv:1507.05464}

\bibitem{Vu}
Vu D.-V.,
Complex Monge-Amp\`ere equation for measures supported on real submanifolds. {\it Preprint} (2016).
{\tt arXiv:1608.02794 }

  
 \bibitem{Xie} 
 Xie J., Periodic points of birational transformations on projective surfaces.
{\it Duke Math. J.}  {\bf 164} (2015), no. 5, 903-932. 

\bibitem{Yomdin}
Yomdin Y., Volume growth and entropy. {\it  Israel J. Math.} {\bf  57}  (1987),  no. 3, 285-300.

\bibitem{Zhang}
Zhang D.-Q.
Invariant hypersurfaces of endomorphisms of projective varieties. 
{\it Adv. Math.} {\bf 252} (2014), 185-203. 

\end{thebibliography}
\end{document}